\documentclass[twoside, 11pt]{article}
\usepackage{mathrsfs,amsfonts,amsmath}
\RequirePackage{ifthen,calc}
\usepackage{theorem}
\usepackage[dvips]{color}

 \renewcommand{\baselinestretch}{1.2}

 \catcode`@=11
 \def\@evenhead{\hbox to\textwidth{\footnotesize\rm\thepage \hfill
  {\it Yuqiang LI}}} 

 \def\@oddhead{\hbox to \textwidth{\footnotesize{\it
   The effects of dynamic branching laws} \hfill\thepage}}

 \renewcommand{\section}{\makeatletter
 \renewcommand{\@seccntformat}[1]{{\csname the##1\endcsname.}\hspace{0.45em}}
 \makeatother \@startsection
{section}
{1}
{0pt}
{\baselineskip}
{0.5\baselineskip}
{\normalsize\bfseries\mathversion{bold}}}

\catcode`@=12

\newtheorem{thm}{\noindent Theorem}[section]
\newtheorem{lem}{\noindent Lemma}[section]
\newtheorem{cor}{\noindent Corollary}[section]
\newtheorem{prop}{\noindent Proposition}[section]

\theoremheaderfont{\normalfont\bfseries} \theorembodyfont{\slshape}
{\theorembodyfont{\rmfamily}} {\theorembodyfont{\rmfamily}
}
{\theorembodyfont{\rmfamily} } {\theorembodyfont{\rmfamily}
}
{\theorembodyfont{\rmfamily} } {\theorembodyfont{\rmfamily}
\newtheorem{rem}{\noindent Remark}[section]}
{\theorembodyfont{\rmfamily} } {\theorembodyfont{\rmfamily} } {\theorembodyfont{\rmfamily}
}

 \def\beqlb{\begin{eqnarray}}\def\eeqlb{\end{eqnarray}}
 \def\beqnn{\begin{eqnarray*}}\def\eeqnn{\end{eqnarray*}}
 \numberwithin{equation}{section}
\def\qed{\hfill$\square$\smallskip}
\def\R{{\mathbb R}}
\def\bfE{{\mathbb{E}}}

\def\e{\mathrm{e}}

\def\d{\mathrm{d}}

\begin{document}

\title{\LARGE\bf Functional ergodic limits of site-dependent branching Brownian motions in $\R$}
\author{ Yuqiang Li \renewcommand{\thefootnote}{\fnsymbol{footnote}}\footnotemark[1]
\\ \small School of Finance and Statistics, East China Normal University,
\\ \small Shanghai 200241, P. R. China.
}
\date{}

\maketitle
\renewcommand{\thefootnote}{\fnsymbol{footnote}}\footnotetext[1]
{Research supported partly by NSFC grant (10901054).}

\centerline{\textbf{Abstract}}

In this paper, we studied the functional ergodic limits of  the site-dependent branching Brownian motions in $\R$.
The results show that the limiting processes are non-degenerate if and only if the variance functions of branching laws are integrable.
When the functions are integrable, although the limiting processes will vary according to the integrals, they are always
positive, infinitely divisible and self-similar, and their marginal distributions are determined by a kind of $1/2$-fractional integral equations.
As a byproduct, the unique non-negative solutions of the integral equations can be explicitly presented by the L\'{e}vy-measure
of the corresponding limiting processes.

\smallskip


\smallskip


\renewcommand{\baselinestretch}{1.2}

\normalsize

\bigskip

\section{Introduction}
By the name ``site-dependent branching Brownian motion" (SDBBM) we mean a  branching particle system where
particles start off at time $t=0$ from a Poisson random field with
Lebesgue intensity measure $\lambda$, move in $\R^d$ according to the Brownian motion, and evolve independently with critical  branching laws at rate $\gamma$. Here the critical branching law of particles at site $x$ is controlled by the generating function
  \beqlb\label{intr-1}
 g(s,x)=s+\sigma(x)(1-s)^2, \qquad 0\leq s\leq1,
 \eeqlb
where $0\leq\sigma(x)\leq 1/2$ is a measurable function on $\R^d$. (\ref{intr-1}) means that a particle at the site $x$ reproduces $0$-offspring with probability $\sigma(x)$, $1$-offspring with probability $1-2\sigma(x)$ and $2$-offsprings with probability $\sigma(x)$. This model generalizes the typical critical branching Brownian motion, which in fact corresponds to the case of $\sigma(x)\equiv 1/2$. It is also a special case of the general branching particle systems discussed in Dynkin \cite{D91}.

Let $N(s)$ denote the random counting measure of a SDBBM at
time $s$, i.e. $N(s)(A)$ is the number of particles in the set
$A$ at time s, and be referred to as the SDBBM for convenience in this paper. Dynkin \cite{D91} had shown that for any non-negative $\phi\in\mathcal{S}(\R^d)$, the space of smooth rapidly decreasing
functions,
 \beqlb\label{intr-2}
 \bfE\Big(\exp\Big\{-\big\langle \phi, N(t)\big\rangle\Big\}\Big)=\exp\Big\{-\int_{\R^d} v(x, t)\d x\Big\},
 \eeqlb
where $0\leq v(x, t)\leq 1$ satisfies the partial differential equation
 \beqlb\label{intr-3}
 \frac{\partial{v}}{\partial t}=\frac{1}{2}\Delta v-\gamma\sigma v^2,
 \eeqlb
 with $v(x, 0)=1-\e^{-\phi(x)}$. Here $\Delta$ denotes the Laplace operator and $\langle f, \mu\rangle=\int f \d \mu$ for any integrable function $f$ on the Borel measure $\mu$. Imagine a SDBBM varied in the way that particles' life become more and more short and the density of initial distribution become more and more high. As a result, the so-called Dawson-Watanabe super-process, say $Y$, appears.  According to Dynkin \cite[Theorem 1.1]{D91}, $Y$ satisfies that for any meaningful initial measure-value $\mu$,
   \beqlb\label{intr-4}
 \bfE_{\mu}(\exp\{-\langle \phi, Y(t)\rangle\})=\exp\Big\{-\int_{\R^d} v(x, t)\mu(\d x)\Big\},
 \eeqlb
where $0\leq v(x, t)\leq 1$ is the solution of (\ref{intr-3})
with $v(x, 0)=\phi(x)$.
From Pinsky \cite[Theorem 4 and Remark 5]{P01} we know that when $d\geq 2$, the invariant measure for $Y$ depends on $\sigma(\cdot)$. However, when $d=1$, the unique invariant measure is $0$, the measure concentrated on the 0-measure, if and only if $\sigma(\cdot)\not=0$ a.s.
It is not surprising that these conclusions also hold for the SDBBMs.

Due to the close relation between the invariant measures and the ergodic limits, naturally, an interesting question arises. How the function $\sigma$ affects the ergodic limits of the SDBBMs or the corresponding super-processes when $d=1$? In this paper, we only discuss the SDBBMs. We believe the same discussions can be moved to the corresponding super-processes.

Let $N$ be a SDBBM. Consider the measure-valued processes
\beqlb\label{s1-2}
X(t)=\int_0^{t}N(s)\d s,\quad t\geq0,
\eeqlb
which is generally referred to as the occupation time (process). Throughout this paper, we understand (\ref{s1-2}) as
 $$\langle X(t), \phi\rangle=\int_0^{t}\langle N(s), \phi\rangle\d s,\quad t\geq0,$$ for any
$\phi\in\mathcal{S}(\R^d)$. The so-called ergodic limits of the SDBBMs are, in general, referred to the limits of $X(T)/T$ as $T\to\infty$.

This study on typical critical branching Brownian motions can be originated to the 70s of the last century. Sawyer and Fleischman \cite{SF78} proved that when $d=1$
 $$
  \lim_{T\to\infty} X(T)(A)<\infty,\;\;\qquad a.s.,
 $$
which immediately leads to
 \beqlb\label{s1-3}
  \lim_{t\to\infty} \frac{X(t)(A)}{t}&=&0,
 \eeqlb
almost surely. Iscoe \cite{I86} showed that  the corresponding super-Brownian motions have the same property when $d=1$. Moreover, Cox and Geriffeath \cite{CG85} found that the same picture arises in voter models.


To step forward, our interests in this paper are to study the  functional ergodic limits of SDBBMs. More precisely, we study the functional limit of $X_T=\{X_T(t),
t\geq0\}$  in $C([0, 1],\mathcal{S}'(\R))$, where
 \beqlb\label{s2-4}
 X_T(t)=\frac{1}{T}\int_0^{Tt} N(s)\d s,
 \eeqlb
and $\mathcal{S}'(\R)$ is the dual space of $\mathcal{S}(\R)$. Our results show that the limiting process is non-degenerate if and only if the variance function is integrable. If $\int_{\R}\sigma(x)\d x=\infty$, then (\ref{s1-3}) holds in probability, and hence,
$X_T$ converges in finite-dimensional distributions to the measure concentrated on the $0$-measure. This result extends the aforementioned result on the typical branching Brownian motions in $\R$ and corresponding super-processes. If $\int_{\R}\sigma(x)\d x<\infty$, then (\ref{s1-3}) is not true. In fact, we prove that the limiting process will vary according to the integral of $\sigma(x)$. But it is always positive, infinitely divisible and self-similar and its marginal distributions are determined by a kind of $1/2$-fractional integral equations. These results are similar in appearance to but different in essence from those in literature (see, for example, \cite{CG85, I86b, T08}) on the typical branching Brownian motions in $\R^2$ and corresponding super-processes.

The methods of this paper consist of two key-points. One is the convergence of finite-dimensional distributions of $X_T$ under the condition of $\int_{\R}\sigma(x)\d x<\infty$. To solve this problem, we need study the convergence of the solutions of some nonlinear integral equations by means of the Gronwall's inequality. Though this idea is similar  to those used in Iscoe \cite{I86, I86b} and Talarczyk \cite{T08}, some
nontrivial modifications are needed to handle the new technical complexities caused
by the site-dependence. The other is to analyze solutions of the $1/2$-fractional integral equations. Based on these analytical results we find the limiting processes are degenerate under the condition of $\int_{\R}\sigma(x)\d x=\infty$. Furthermore, by using the Levy-Khintchine representation, we show the positivity of the limiting process under the condition of $\int_{\R}\sigma(x)\d x<\infty$ and get the explicit expressions of the solutions of the integral equations. This trick differs from that used by Iscoe \cite{I86b}.

There is much literature on occupation times of branching particle systems and the closely related super-processes; see, for example, Dawson {\it et al} \cite{DG01} and the references therein. Studying the functional limits of occupation times of typical $(d,\alpha,\beta)$-branching particle systems was triggered by Bojdecki {\it et al} \cite{BGT061} and was developed and generalized by, for example, Bojdecki {\it et al} \cite{BGT071}, Li \cite{L10}, Li and Xiao \cite{LX09} and Milo\'{s} \cite{M07}. Li \cite{Li11a, Li11c} introduced a kind of site-dependent branching particle systems which are same as the SDBBMs except that the particles move according to the $\vec{\alpha}=(\alpha_1,\cdots,\alpha_d)$-stable L\'{e}vy motion whose $i$-th component, $i=1,2,\cdots, d$, is an $\alpha_i$-stable Levy process and independent of other components. Under the assumption that $0<\int_{\R}\sigma(x)\d x<\infty$, the author has studied the functional theorems of central limit type of the occupation times except the case of $\sum_{i=1}^d1/\alpha_i=1/2$ and obtained some interesting results which differ from the existing results of the typical $(d, \alpha,\beta)$-branching systems and the particle systems without branching. Observe that the SDBBM in $R$ is in essence the site-dependent branching particle system in the case of $\sum_{i=1}^d1/\alpha_i=1/2$. This paper also completes the picture of research in this direction.

Without other statement, in this paper, we use $M$ to denote an unspecified positive finite constant which may not necessarily be
the same in each occurrence.

The remainder of this paper is organized as follows. In Section 2,
we report the main results. Section 3 devotes to studying the finite-dimensional distribution of $X_T$ and other related lemmas. Sections 4 and 5 include the proofs of Theorem \ref{thm1} and Theorem \ref{thm5}, respectively. In the last section, i.e. Section 6, the limiting process $\xi$ is discussed.

\section{Main results}
Consider a SDBBM $N=\{N(s)\}$ in $\R$. Let $B=\{B(s)\}$ denote the corresponding Brownian motion.
$\{L_t\}_{t\geq0}$ denotes the semi-group of Brownian motion. Then
 \beqlb\label{s2-0}
 L_sf(x):=\bfE(f(\xi(t+s))|\xi(t)=x)=\int_{\R}\frac{1}{\sqrt{2\pi s}}\exp\Big\{-\frac{(y-x)^2}{2s}\Big\}f(y)\d y,
 \eeqlb
for all $s, t\geq0$, $x\in\R$ and bounded measurable functions
$f$. To avoid misunderstanding, we sometimes write $L_sf(x)$ as
$L_s(f(\cdot))(x)$. Let
 $$p_s(x):=\frac{1}{\sqrt{2\pi s}}\exp\Big\{-\frac{x^2}{2s}\Big\},$$  for all $s>0$, $x\in\R$.  Recall that
 \beqlb\label{s2-0-1}
 p_{rt}(x)=r^{-1/2}p_t(r^{-1/2}x).
 \eeqlb
Therefore,
 \beqlb\label{s2-0-2}
 L_{rt}f(r^{1/2}x)=\int_{\R}r^{-1/2}p_t(x-r^{-1/2}y)f(y)\d y,
 \eeqlb
for all $r>0$. Now we define the rescaled occupation time process $X_T=\{X_T(t),
t\geq0\}$ as (\ref{s2-4}). The main results of this paper read as follows.

\begin{thm}\label{thm4}
If $\int_{\R}\sigma(x)\d x=\infty$, then $X_T$ converges to $0$ in finite-dimensional distributions as $T\to\infty$.
\end{thm}

\begin{thm}\label{thm1} Suppose $K=\gamma\int_{\R}\sigma(x)\d x<\infty$. Then as $T\to\infty$, the rescaled occupation time process $\{X_T(t),, 0\leq t\leq 1\}$ converges weakly in $C([0, 1],\mathcal{S}'(\R))$ to a process $X=\lambda\xi$ where $\xi=\{\xi(t)\}_{t\in[0, 1]}$ is a non-negative process whose Laplace transforms of finite-dimensional distributions
 \beqlb\label{thm1-1}
   \bfE\Big(\exp\Big\{-\sum_{k=1}^n \theta_k\xi(t_k)\Big\}\Big)=\exp\bigg\{K\int_{1-t_n}^{1}\Lambda^2(s)\d s-\sum_{k=1}^nt_k\theta_k
\bigg\},
  \eeqlb
for any given non-negative constants $\theta_1,\theta_2,\cdots,\theta_n$ and $0\leq t_1<t_2<\cdots<t_n\leq 1$, where $\Lambda(s)$ is the unique nonnegative solution of the equation
 \beqlb\label{thm1-2}
 y(s)=\sum_{k=1}^n\theta_k\int_0^s\frac{{\bf 1}_{[0, t_k]}(1-u)}{\sqrt{2\pi (s-u)}}\d u-K \int_0^{s}\frac{y^2(u)\d u}{\sqrt{2\pi (s-u)}},\quad s\in[0, 1].
 \eeqlb
\end{thm}

From Theorem \ref{thm1}, it is easy to see that
\begin{prop}\label{prop4}
The process $\xi$ is infinitely divisible, self-similar with index $1$, non-decreasing and nonnegative, and has continuous paths.
\end{prop}

To get more information about $\xi$, we need to make more careful study on the solution of (\ref{thm1-2}). For any given $\theta>0$, let $\Lambda(s,\theta)$, $0\leq s\leq 1,$ be the unique non-negative solution of the equation
 \beqlb\label{thm5-1}
 y(s)=\sqrt{\frac{2s}{\pi}}\theta-K \int_0^{s}\frac{y^2(u)\d u}{\sqrt{2\pi (s-u)}},
 \eeqlb
which is the special case of (\ref{thm1-2}) with $k=1$, $\theta_1=\theta$ and $t_1=1$.
We obtain the following results.
\begin{thm}\label{thm5}
 {\rm (1)} $\Lambda(s,\theta)$ is continuous, differentiable and non-decreasing on $s\in[0, 1]$.
Moreover, $\Lambda(s,\theta)\leq\sqrt{\theta/K}$.

{\rm (2)} $\Lambda(s,\theta)$ can be continuously extended to the unique non-decreasing and non-negative solution of (\ref{thm5-1}) for all $s\geq0$. Denote the extension by $\Lambda(s,\theta)$ as well. Then
 \beqlb\label{thm5-2}
 \Lambda(s,\theta)=\frac{1}{\sqrt{s}}\Lambda(1, \theta s)
 \eeqlb
 and $\lim_{s\to\infty}\Lambda(s, 1)=1/\sqrt{K}$.
\end{thm}
 \begin{prop}\label{prop1}
 There exists a measure $\nu$ on $(0,\infty)$ with $\nu((0,\infty))=\infty$ and $$\int_0^\infty x\nu(\d x)=1,\qquad\int_0^\infty x^2\nu(\d x)=2K/\pi,$$
 such that for any $\theta>0$
 \beqlb\label{prop-1}
 \bfE\big(\e^{-\theta \xi(t)}\big)=\exp\Big\{-\int_0^\infty(1-\e^{-t\theta x})\nu(\d x)\Big\}.
 \eeqlb
Furthermore, for any $s\geq0$ and $\theta>0$,
 \beqlb\label{thm5-4}
 K\Lambda^2(s,\theta)=\theta-\int_0^\infty\theta\e^{-\theta s x}\nu(\d x),
 \eeqlb
and
 \beqlb\label{thm5-3}
 \Lambda(s,\theta)=\frac{1}{\sqrt{2\pi s}}\int_0^\infty Q(\theta s x)\nu(\d x),
 \eeqlb
 where
 $$Q(x)=\frac{\sqrt{x}}{\e^x}\int_0^x\frac{\e^y}{\sqrt{y}}\d y,\;\; x>0.$$
\end{prop}

From Proposition \ref{prop1} we immediately get that
 \begin{cor}
 For any $t>0$, $\xi(t)$ is non-trivial and positive.
 \end{cor}

At the end of this section, let us make some remarks on the results mentioned above.

\begin{rem} (1) Cox and Geriffeath \cite{CG85} proved for the typical critical branching Brownian motions by a method of cumulants that when the spatial dimension $d=2$,
 \beqlb\label{s1-4}
 \frac{X_T(1)(A)}{T\lambda(A)}\to\varsigma, \;\;\; \text{in law},
 \eeqlb
where $\varsigma$ is a nontrivial infinitely
divisible random variable with moments of all order. Via Laplace transforms of measure-valued random variables and nonlinear partial differential equations, Iscoe \cite{I86b} proved the same result for the corresponding super-Brownian motions
and further pointed out the positivity of $\varsigma$ in (\ref{s1-4}). When $t=1$ is fixed, Theorem \ref{thm1} and Proposition \ref{prop1} reveal the similar results. However, the methods used in this paper to prove the positivity are different from those used in Iscoe \cite{I86b}.

(2) The similar functional ergodic theorems of the typical branching Brownian motion and the related super-Brownian motions were reported by Talarczyk \cite{T08} and Iscoe \cite{I86b}, respectively, in the case of the spatial dimension $d=2$. Compared with those results, our result is essentially different in the Laplace functions of the limit processes. In fact, using the Riemann-Liouville type $1/2$-fractional integral operator,
 $$J^{1/2}(f)(s):=\int_0^{s}\frac{f(u)\d u}{\sqrt{2\pi (s-u)}},$$
the equation (\ref{thm1-2}) can be rewritten as
 \beqnn
 \Lambda(s)+KJ^{1/2}\Lambda^2(s)=\sum_{k=1}^n\theta_k\int_0^s\frac{{\bf 1}_{[0, t_k]}(1-s+u)}{\sqrt{2\pi u}}\d u=J^{1/2}h(s),
 \eeqnn
where
 $$h(s)=\sum_{k=1}^n\theta_k{\bf 1}_{[0, t_k]}(1-s).$$
This equation can not be inferred from the corresponding results in Talarczyk \cite{T08} and Iscoe \cite{I86b}. In addition, (\ref{thm5-4}) and (\ref{thm5-3}) show the one-on-one corresponding relations between the non-negative solution of (\ref{thm5-1}) and the L\'{e}vy measure of $\xi$.

(3) By some basic renewal discussion, it is easy to obtain that
  \beqlb\label{s2-3}
 \bfE(\langle N(s), \phi\rangle)=\int_{\R} L_s\phi(x)\d x=\langle \lambda, \phi\rangle.
 \eeqlb
To assure non-degenerate limits existing, the  rescaled occupation time fluctuation processes should be defined as follows.
 \beqlb\label{s2-5}
 \langle Y_T(t),\phi\rangle=\frac{1}{F_T}\int_0^{Tt}\langle N(s)-\lambda,\phi\rangle,
 \eeqlb
where $$F_T=\begin{cases} T ,\qquad &0<\int_{\R}\sigma(x)\d x<\infty;
\\ T^{3/4}, \qquad &\int_{\R}\sigma(x)\d x=0.\end{cases}$$
For the rescaled occupation time fluctuations, we have the following functional limits.
\begin{itemize}
\item[(i)] When $0<\int_{\R}\phi(y)\d y<\infty$, the rescaled occupation time fluctuation process $\{Y_T(t),, 0\leq t\leq 1\}$ converges weakly in $C([0, 1],\mathcal{S}'(\R))$ to the process $Y=\lambda\eta$ where $\eta(t)=\xi(t)-t$ and $\xi$ is the process in Theorem \ref{thm1}.
\item[(ii)] When $\int_{\R}\phi(y)\d y=0$, i.e. $\phi(y)=0$ a.e., the rescaled occupation time fluctuation process $\{Y_T(t),, 0\leq t\leq 1\}$ converges weakly in $C([0, 1],\mathcal{S}'(\R))$ to the process $Y=k\lambda\eta$ where $\eta=\{\eta(t)\}_{t\in[0, 1]}$ is the fractional Brownian motion with Hurst index $\frac{3}{4}$ and $k$ is a constant.
\end{itemize}
The part (i) is an immediate result from Theorem \ref{thm1}, and the part (ii) was essentially investigated in Bojdecki {\it et al} \cite{BGT061}.
\end{rem}

\bigskip

\section{Finite-dimensional distributions of $X_T$}
Define a sequence of random variables
$\tilde{X}_T$ in $\mathcal{S}'(\R^{2})$ as follows: For any
$n\geq0$, let
 \beqlb\label{s3-0}
 \langle\tilde{X}_T, \psi\rangle=:\int_0^1\langle X_T(t), \psi(\cdot, t)\rangle\d t=\int_0^T\langle N(s), \psi_T(\cdot, s)\rangle\d s,
 \eeqlb
where $\psi\in\mathcal{S}(\R^{2})$ and
 \beqlb\label{s3-1}
 \psi_T(x, s)=\frac{1}{T}\int_{s/T}^1\psi(x, t)\d t, \qquad\text{for}\;\; s\in[0, T].
 \eeqlb
\begin{lem}\label{s3-lem1} For any nonnegative $\psi\in\mathcal{S}(\R^2)$,
 \beqlb\label{lem0-3}
  \bfE\Big(\exp\Big\{-\langle\tilde{X}_T, \psi\rangle\Big\}\Big)&=&\exp\bigg\{-\int_{\R} V_{\psi_T}(x, T, 0)\d x\bigg\},
  \eeqlb
 where $\psi_T(x, s)$ is defined by (\ref{s3-1}) and $V_{\psi_T}(x, t, r)$ is a continuous function defined on $\R\times\{(t, r): t\geq0, r\geq0 , t+r\leq T\}$ and satisfies that
  \beqlb\label{s3-9}
 V_{\psi_T}(x, t, r)&=&\int_0^tL_s\Big(\psi_T(\cdot, r+s)\big(1-V_{\psi_T}(\cdot, t-s, r+s)\big)\Big)(x)\d s\nonumber
 \\&&-\gamma \int_0^tL_s\Big(\sigma(\cdot)V_{\psi_T}^2(\cdot, t-s,
 r+s)\Big)(x)\d s.\qquad
 \eeqlb
\end{lem}
{\bf Proof}. On $\R\times\{(t, r): t\geq0, r\geq0 , t+r\leq T\}$, define
 \beqlb\label{s3-2}
 H_{\psi_T}(x, t, r):=\bfE_x\Big(\exp\big\{-\int_0^t\big\langle N(s), \psi_T(\cdot, r+s)\big\rangle\d s\big\}\Big).
 \eeqlb
Since $N_0$ is a Poisson random measure with Lebesgue intensity
measure, from (\ref{s3-0}), it follows that
 \beqlb\label{s3-3}
 \bfE\Big(\e^{-\langle \tilde{X}_T, \psi\rangle}\Big)=\exp\Bigg\{\int_{\R}\big[H_{\psi_T}(x,T,0)-1\big]\d x
 \Bigg\}.\qquad
 \eeqlb
By renewal
arguments, (\ref{s3-2}) implies that
 \beqlb\label{s3-7}
 H_{\psi_T}(x, t, r)&=&\e^{-\gamma t}\bfE_x\Bigg\{\exp\bigg(-\int_0^t\psi_T(\xi(s), r+s)\d
 s\bigg)\Bigg\}\nonumber
 \\&&+\int_0^t\gamma\e^{-\gamma s}\bfE_x\Bigg\{\exp\bigg(-\int_0^s\psi_T(\xi(u), r+u)\d u\bigg)\nonumber
 \\&&\times\Bigg[\bfE_{\xi(s)}\exp\bigg(-\int_0^{t-s}\big\langle N(u),\psi_T(\cdot, r+s+u)\big\rangle\d u\bigg)\Bigg]
 ^{k(\xi(s))}\Bigg\}\d s,\qquad
 \eeqlb
where $k(\xi(s))$ denotes the number of particles generated at the
first splitting time. Note that the process $k(x)$ is independent of the
Brownian motion $\xi$ and for any $0<z<1$
 \beqlb\label{s3-7-1}
 \bfE(z^{k(x)})=g(z, x)=z+\sigma(x)(1-z)^2.
 \eeqlb
(\ref{s3-7}) yields that
 \beqlb\label{s3-8}
 &&H_{\psi_T}(x, t, r)=\e^{-\gamma t}I_{\psi_T}(x, t, r)+
 \int_0^t\gamma\e^{-\gamma(t-s)}K_{\psi_T}(x, t-s, r, s)\d s,
 \eeqlb
where
 \beqnn
 I_{\psi_T}(x, t, r)&=&\bfE_x \exp\bigg(-\int_0^t\psi_T(\xi(s), r+s)\d
 s\bigg),
  \\K_{\psi_T}(x, t, r, s)&=&\bfE_x \bigg[\exp\bigg(-\int_0^t
  \psi_T(\xi(u), r+u)\d u\bigg)g\big(H_{\psi_T}(\xi(t),s,r+t),
 \xi(t)\big)\bigg].
  \eeqnn
for any $(x, t, r)\in \R\times\{(t, r): t\geq0, r\geq0 , t+r\leq T\}$.  By the Feynman-Kac formula,
 \beqnn
 \frac{\partial I_{\psi_T}}{\partial t}&=&\Big(\frac{1}{2}\Delta+\frac{\partial}{\partial r}-\psi_T(x,
 r)\Big) I_{\psi_T},
 \\ \frac{\partial K_{\psi_T}}{\partial t}&=&\Big(\frac{1}{2}\Delta+\frac{\partial}{\partial r}-\psi_T(x,
 r)\Big)K_{\psi_T}.
 \eeqnn
Therefore, (\ref{s3-8}) indicates that
 \beqnn
 \frac{\partial
H_{\psi_T}(x, t, r)}{\partial t}&=&\Big(\frac{1}{2}\Delta+\frac{\partial}{\partial r}-\psi_T(x,
 r)\Big)H_{\psi_T}(x, t, r)+\gamma\bfE_x\Big[g\big(H_{\psi_T}(\xi(0),t,r),
 x\big)\Big],\nonumber
 \eeqnn
which plus (\ref{s3-7-1}) leads to that
 \beqlb\label{lem0-2}
 \frac{\partial
H_{\psi_T}(x, t, r)}{\partial
t}&=&\Big(\frac{1}{2}\Delta+\frac{\partial}{\partial r}-\psi_T(x,
 r)\Big)H_{\psi_T}(x, t, r)+
 \gamma\sigma(x)(1-H_{\psi_T}(x,t,r))^2.\qquad
 \eeqlb
Let
 \beqlb\label{s3-4}
 V_{\psi_T}(x, t, r)=1-H_{\psi_T}(x, t, r).
 \eeqlb
Then it is easy to see that $V_{\psi_T}(x, t, r)$ is continuous on $\R\times\{(t, r): t\geq0, r\geq0 , t+r\leq T\}$. Furthermore, from (\ref{lem0-2}) it follows that
 \beqlb\label{s3-5}
\frac{\partial V_{\psi_T}(x, t, r)}{\partial t}&=&
 \Big(\frac{1}{2}\Delta+\frac{\partial}{\partial r}\Big)V_{\psi_T}(x, t, r)-\gamma\sigma(x)V_{\psi_T}^2(x, t, r)\nonumber
 \\&&+\psi_T(x, r)\big(1-V_{\psi_T}(x, t, r)\big),
 \eeqlb
which implies (\ref{s3-9}). Moreover, from (\ref{s3-3}) and (\ref{s3-4}), (\ref{lem0-3}) follows.  \qed

 Since $X_T$ is a $\mathcal{S}'(\R)$-valued process, as is well-known, the finite dimensional distributions of $X_T$ is determined by the family of Laplace functions, i.e.,
 \beqnn
 \bfE\Big(\exp\big\{-\sum_{k=1}^n\langle X_T(t_k), \phi_k\rangle\big\}\Big),
 \eeqnn
for any given nonnegative $\phi_1,\phi_2,\cdots,\phi_n\in\mathcal{S}(\R)$, $0<t_1<t_2<\cdots<t_n\leq 1$ and $n\geq 1$.
 \begin{lem}\label{s3-lem3}
  For $0\leq t_1<t_2<\cdots<t_n\leq 1$, nonnegative $\phi_1,\phi_2,\cdots,\phi_n\in\mathcal{S}(\R)$, and $T>0$,
  \beqlb\label{s3-lem3-1}
  \bfE\Big(\exp\Big\{-\sum_{k=1}^n\langle X_T(t_k), \phi_k\rangle\Big\}\Big)&=&\exp\bigg\{-\int_{\R} V_{\psi_T}(x, T, 0)\d x\bigg\},
  \eeqlb
 where $0\leq  V_{\psi_T}(x, T, 0)\leq 1$ satisfies the equation (\ref{s3-9}) with
  \beqlb\label{s3-lem3-2}
  \psi_T(x, t)=\frac{1}{T}\psi(x, \frac{t}{T}),\qquad \psi(x, t)=\sum_{k=1}^n\phi_k(x){\bf 1}_{[0, t_k]}(t).
  \eeqlb
  \end{lem}
 {\bf Proof}\;\; There exists $\psi_{T,m}$ of the form (\ref{s3-1}) such that for any $(x, t)$, $\psi_{T,m}(x, t)$ converges to $\psi_T(x, t)$ in the monotone decreasing way  as $m\to\infty$. Lemma \ref{s3-lem3} follows from Lemma \ref{s3-lem1}. The details are same as those lead to Lemma 2.5 in \cite{T08} and omitted.\qed

\begin{rem}\label{s3-rem}
Let $v(x, t)=V_{\psi_T}(x, t, T-t)$. Lemma \ref{s3-lem3} and (\ref{s3-5}) imply that $v(x, t)$ satisfies that
  \beqnn\label{s3-rem-1}
\frac{\partial v(x, t)}{\partial t}&=&
 \frac{1}{2}\Delta v(x, t)-\gamma\sigma(x)v^2(x, t)+\frac{1}{T}\sum_{k=1}^n\phi_k(x){\bf 1}_{[0, t_k]}(\frac{T-t}{T})\big(1-v(x, t)\big),\quad
 \eeqnn
with $v(x,0)=0$, and that
 \beqnn\label{s3-rem-2}
  \bfE\Big(\exp\Big\{-\sum_{k=1}^n\langle X_T(t_k), \phi_k\rangle\Big\}\Big)&=&\exp\bigg\{-\int_{\R} v(x, T)\d x\bigg\}.
  \eeqnn
\end{rem}

\section{Weak convergence of $X_T$ when $\int_{\R}\sigma(x)\d x<\infty$.}
The aim of this section is to prove the weak convergence of the rescaled occupation time $X_T$. Note that the weak convergence follows from the convergence of finite-dimensional distributions plus tightness.

\begin{thm}\label{thm2}
 As $T\to\infty$, the occupation time process $\{X_T(t),, 0\leq t\leq 1\}$ converges in finite-dimensional distributions to a $\mathcal{S}'(R)$-valued process $X=\lambda\xi$ where $\xi=\{\xi(t)\}_{t\in[0, 1]}$ is a non-negative process whose Laplace transforms of finite-dimensional distributions are as follows.
 \beqlb\label{thm2-1}
   \bfE\Big(\exp\Big\{-\sum_{k=1}^n \theta_k\xi(t_k)\Big\}\Big)=\exp\bigg\{K\int_0^{1}\Lambda^2(s)\d s-\sum_{k=1}^nt_k\theta_k
\bigg\},
  \eeqlb
for any given non-negative constants $\theta_1,\theta_2,\cdots,\theta_n$ and $0\leq t_1<t_2<\cdots<t_n\leq 1$, where $\Lambda(s)$ is the unique nonnegative solution of the equation (\ref{thm1-2}).
\end{thm}

To prove this result, we need some auxiliary lemmas. For this end, we first remark that, throughout this section, we always assume that $\phi_1,\phi_2,\cdots,\phi_n\in\mathcal{S}(\R)$ are arbitrary non-negative functions, $0\leq t_1<t_2<\cdots<t_n\leq 1$, and $\psi(x, t)$, $\psi_{T}(x, t)$ have the forms in (\ref{s3-lem3-2}) and $K=\gamma\int_{\R}\sigma(x)\d x<\infty$. Define
 \beqlb\label{s4-4}
  \bar{H}_T(x, s)&:=&T^{1/2}V_{\psi_T}(x, Ts, T(1-s)).
 \eeqlb
Observing (\ref{s3-9}), we derive that
 \beqlb
\bar{H}_T(x, s)&=&T^{1/2}\int_0^{Ts}L_u\Big(\psi_T(\cdot, T(1-s)+u)\nonumber
 \big(1-V_{\psi_T}(\cdot, Ts-u, T(1-s)+u)\big)\Big)(x)\d u\nonumber
 \\&&-\gamma T^{1/2}\int_0^{Ts}L_u\Big(\sigma(\cdot)V_{\psi_T}^2(\cdot, Ts-u, T(1-s)+u)\Big)(x)\d u\nonumber
 \\&=&T^{3/2}\int_0^{s}L_{Tu}\psi_T(x, T(1-s+u))\d u-\gamma T^{1/2}\int_0^{s}L_{Tu}(\sigma(\cdot)\bar{H}_T^2(\cdot, s-u))(x)\d u\nonumber
 \\&&-\int_0^{s}TL_{Tu}(\psi_T(\cdot, T(1-s+u))\bar{H}_T(\cdot, s-u))(x)\d u.\nonumber
  \eeqlb
Substituting (\ref{s3-lem3-2}) into the above formula, we further have that
 \beqlb\label{s4-5}
 \bar{H}_T(x, s)&=&I_1(T, x, s)-I_2(T, x, s)-I_3(T, x, s),
 \eeqlb
where
 \beqlb
 I_1(T, x, s)&=&T^{1/2}\int_0^{s}L_{Tu}\psi(x, 1-s+u )\d u\nonumber
 \\&=&\int_0^{s}\d u\int_{\R}\frac{1}{\sqrt{2\pi u}}\e^{-\frac{(x-y)^2}{2Tu}}\psi(y, 1-s+u)\d y;\label{s3-2-1}
 \\I_2(T, x, s)&=&\int_0^{s}L_{Tu}(\psi(\cdot, 1-s+u)\bar{H}_T(\cdot, s-u))(x)\d u;\label{s3-2-2}
 \\I_3(T, x, s)&=&\gamma T^{1/2}\int_0^{s}L_{Tu}(\sigma(\cdot)\bar{H}_T^2(\cdot, s-u))(x)\d u.\label{s3-2-3}
 \eeqlb
 \begin{lem}\label{lem7}
There exists a constant $M$ such that for any $T>0$ and $(x, s)\in\R\times[0, 1]$,
 $$
\bar{H}_T(x, s)\leq M.
 $$
\end{lem}
{\bf Proof}. From (\ref{s2-3}), (\ref{s3-4}) and (\ref{s3-lem3-2}) it follows that
 \beqnn
 \bar{H}_T(x, s)&=&T^{1/2}V_{\psi_T}(x, Ts, T(1-s))\nonumber
 \\&\leq& T^{1/2}\int_0^{Ts}L_u\psi_{T}(x, T(1-s)+u)\d u\nonumber
 \\&=& T^{1/2}\int_0^{s}L_{Tu}\psi(x, 1-s+u)\d u.
 \eeqnn
By using (\ref{s2-0}), we further have that
 \beqnn
 \bar{H}_T(x, s)&\leq& T^{1/2}\int_0^{s}\d u\int_{\R}\frac{1}{\sqrt{2\pi Tu}}\e^{-\frac{(x-y)^2}{2Tu}}\psi(y, 1-s+u)\d y
 \\&\leq & \int_0^{s}\frac{1}{\sqrt{2\pi u}}\d u\int_{\R}\psi(y, 1-s+u)\d y.
 \eeqnn
 Recall that $\psi=\sum_{k=1}^n\phi_k(x){\bf 1}_{[0, t_k]}(t)$ and $\phi_k\in\mathcal{S}(\R)$. There exists a positive constant $M$, such that for any $T>0$ and $(x, s)\in\R\times[0, 1]$,
 \beqlb\label{lem7-1}
\bar{H}_T(x, s)\leq M\int_0^s\frac{1}{2\sqrt{u}}\d u=M\sqrt{s}\leq M.
 \eeqlb

\begin{lem}\label{lem1}
For any $(x, s)\in\R\times[0, 1]$, as $T\to\infty$
 \beqlb\label{lem1-1}
  I_1(T, x, s)\to\sum_{k=1}^n\int_{\R}\phi_k(y)\d y\int_0^s\frac{{\bf 1}_{[0, t_k]}(1-u)}{\sqrt{2\pi (s-u)}}\d u.
 \eeqlb
\end{lem}
{\bf Proof}\;\; From (\ref{s3-2-1}), it is easy to see that for all $T>0$ and $(x, s)\in\R\times[0, 1]$,
 $$
 I_1(T, x, s)=\sum_{k=1}^n\int_0^{s}\frac{{\bf 1}_{[0, t_k]}(1-s+u)}{\sqrt{2\pi u}}\d u\int_{\R}\e^{-\frac{(x-y)^2}{2Tu}}\phi_k(y)\d y.
 $$
Letting $T\to\infty$, by the dominated convergence theorem, we can readily get 
(\ref{lem1-1}). \qed

\begin{lem}\label{lem2}
For all $T>0$ and $(x, s)\in\R\times[0, 1]$, there exists a constant $M>0$ such that as $T\to\infty$,
 \beqlb\label{lem2-1}
  M\sum_{k=1}^n\int_{\R}\phi_k(x)\d x\int_0^{s}\frac{\d u}{\sqrt{2\pi Tu}}\geq I_2(T, x, s)\to 0.
 \eeqlb
 \end{lem}
{\bf Proof}\;\; From (\ref{s3-2-2}) and (\ref{lem7-1}) we obtain that
 $$I_2(T, x, s)\leq M\int_0^{s}L_{Tu}\psi(x, 1-s+u)\d u,$$
which and (\ref{s2-0}) further imply that
  $$I_2(T, x, s)\leq M\int_0^{s}\frac{1}{\sqrt{2\pi Tu}}\d u\int_{\R}\psi(y, 1-s+u)\d y.$$
Now substituting (\ref{s3-lem3-2}) into the above formula immediately leads to (\ref{lem2-1}).\qed

\begin{lem}\label{lem3} There exists  a bounded and measurable function $G(s)$ such that for any $s\in[0, 1]$, as $T\to\infty$
 \beqlb\label{lem3-1}
 \int_{\R}\sigma(x)\bar{H}_T^2(x, s)\d x\to G(s).
 \eeqlb
\end{lem}
{\bf Proof}\;\;To simplify the notation, let
 \beqlb\label{lem3-0}
 G_T(x, s)=\sigma(x)\bar{H}^2_T(x, s).
 \eeqlb
For any given $s\in[0, 1]$, define the distance between $G_{T_1}(x, s)$ and $G_{T_2}(x, s)$ for any $0<T_1<T_2$ as follows.
 \beqnn\label{lem3-2}
 d_{G}(T_1, T_2; s)&=&\int_{\R}|G_{T_1}(x, s)-G_{T_2}(x, s)|\d x\nonumber
 \\&=&\int_{\R}\sigma(x)\Big|\bar{H}_{T_1}(x, s)-\bar{H}_{T_2}(x, s)\Big|\Big|\bar{H}_{T_1}(x, s)+\bar{H}_{T_2}(x, s)\Big|\d x.
 \eeqnn
By using (\ref{lem7-1}), we know that
 \beqlb\label{lem3-3}
  d_{G}(T_1, T_2; s)\leq M\int_{\R}\sigma(x)\Big|\bar{H}_{T_1}(x, s)-\bar{H}_{T_2}(x, s)\Big|\d x.
 \eeqlb
Substituting (\ref{s4-5}) into the right hand side of (\ref{lem3-3}), we further obtain that
 \beqlb\label{lem3-4}
  d_{G}(T_1, T_2; s)&\leq& \bar{\rho}_{G}(T_1, T_2; s)+M\int_{\R}\sigma(x)\Big|I_3(T_1, x, s)-I_3(T_2, x, s)\Big|\d x,
 \eeqlb
where
 \beqlb\label{lem3-5}
 \bar{\rho}_{G}(T_1, T_2; s)=M\int_{\R}\sigma(x)\sum_{i=1}^2\Big|I_i(T_1, x, s)-I_i(T_2, x, s)\Big|\d x.
 \eeqlb
Furthermore, from (\ref{s3-2-3}) and (\ref{lem3-0})  it follows that
 \beqlb\label{lem3-6}
 &&M\int_{\R}\sigma(x)\Big|I_3(T_1, x, s)-I_3(T_2, x, s)\Big|\d x\leq \tilde{\rho}_{G}(T_1,T_2; s)+\hat{\rho}_G(T_1,T_2; s),\nonumber
 \eeqlb
 where
  \beqlb\label{lem3-7}
 \tilde{\rho}_{G}(T_1,T_2; s)&=&M\gamma \int_{\R}\sigma(x)\d x\int_0^{s}\Big|T_1^{1/2}L_{T_1 u}G_{T_2}(x, s-u)-T_2^{1/2}L_{T_2 u}G_{T_2}(x, s-u)\Big|\d u\nonumber
 \\&=&M \int_{\R}\sigma(x)\d x\int_0^{s}\Big|\int_{\R}\frac{1}{\sqrt{2\pi u}}\e^{-\frac{(y-x)^2}{2T_2u}}G_{T_2}(y, s-u)\d y\nonumber
 \\&&\qquad\qquad\qquad\qquad-\int_{\R}\frac{1}{\sqrt{2\pi u}}\e^{-\frac{(y-x)^2}{2T_1u}}G_{T_2}(y, s-u)\d y\Big|\d u\nonumber
 \\&=& M \int_{\R}\sigma(x)\d x\int_0^{s}\frac{\d u}{\sqrt{2\pi u}}\int_{\R}\Big(\e^{-\frac{(y-x)^2}{2T_2u}}-\e^{-\frac{(y-x)^2}{2T_1u}}\Big)G_{T_2}(y, s-u)\d y,
 \eeqlb
and
 \beqlb\label{lem3-8}
  \hat{\rho}_G(T_1,T_2; s)&=&M\gamma T_1^{1/2}\int_{\R}\sigma(x)\d x\int_0^{s}L_{T_1 u}\Big|G_{T_1}(\cdot, s-u)-G_{T_2}(\cdot, s-u)\Big|(x)\d u\nonumber
  \\&\leq&M\int_0^{s}\frac{\d u}{\sqrt{2\pi u}}\int_{\R}\Big|G_{T_1}(y, s-u)-G_{T_2}(y, s-u)\Big|\d y\nonumber
  \\&\leq&M\int_0^{s}\frac{d_{G}(T_1, T_2; u)}{\sqrt{2\pi (s-u)}}\d u.
 \eeqlb
Furthermore, substituting (\ref{lem3-0}) and (\ref{lem7-1}) into (\ref{lem3-7}) yields that
 \beqlb\label{lem3-9}
 \tilde{\rho}_{G}(T_1,T_2; s)&\leq& M \int_{\R}\sigma(x)\d x\int_0^{s}\frac{\d u}{\sqrt{2\pi u}}\int_{\R}\Big(1-\e^{-\frac{(y-x)^2}{2T_1u}}\Big)\sigma(y)\d y\nonumber
 \\&\leq&M \int_{\R}\sigma(x)\d x\int_0^{1}\frac{\d u}{\sqrt{2\pi u}}\int_{\R}\Big(1-\e^{-\frac{(y-x)^2}{2T_1u}}\Big)\sigma(y)\d y=:\delta_1(T_1),\qquad
 \eeqlb
for any $T_2>T_1$  and $s\in[0, 1]$. It is easy to see that as $T_1\to\infty$ ,
 \beqlb\label{lem3-10}
 \delta_1(T_1)\to 0.
 \eeqlb
On the other hand, applying (\ref{s3-2-1}) and (\ref{lem2-1}) to (\ref{lem3-5}), we have that
 \beqlb\label{lem3-11}
 \bar{\rho}_{G}(T_1, T_2; s)&\leq&M\int_{\R}\sigma(x)\d x\int_0^s\frac{\d u}{\sqrt{2\pi u}}\int_{\R}[\e^{-\frac{(x-y)^2}{2T_2u}}-\e^{-\frac{(x-y)^2}{2T_1u}}]\psi(y, 1-s+u)\d y\nonumber
 \\&&+M\int_{\R}\sigma(x)(I_2(T_1, x, s)+I_2(T_2, x, s))\d x\nonumber
 \\&\leq& M\int_{\R}\sigma(x)\d x\int_0^{s}\frac{1}{\sqrt{2\pi u}}\d u\sum_{k=1}^n\int_{\R}[1-\e^{-\frac{(x-y)^2}{2T_1u}}]\phi_k(y)\d y\nonumber
 \\&&+M\int_{\R}\sigma(x)\d x\int_0^{s}\frac{\d u}{\sqrt{2\pi T_1u}}\sum_{k=1}^n\int_{\R}\phi_k(y)\d y.\nonumber
 \eeqlb
 Since $\phi_k$, $k=1,2,\cdots,n$,  all are in $\mathcal{S}(\R)$, there exists a constant $M>0$ such that
 \beqlb\label{lem3-12}
 \bar{\rho}_{G}(T_1, T_2; s)&\leq& M\int_{\R}\sigma(x)\d x\int_0^{1}\frac{1}{\sqrt{2\pi u}}\d u\int_{\R}[1-\e^{-\frac{(x-y)^2}{2T_1u}}]\phi(y)\d y\nonumber
 \\&&+M\int_{\R}\sigma(x)\d x\int_{\R}\phi(y)\d y\int_0^{1}\frac{1}{\sqrt{2\pi T_1u}}\d u\nonumber
 \\&=:&\delta_2(T_1)\to 0,
 \eeqlb
as $T_1\to\infty$. Let
\beqlb\label{lem3-14}
\delta(T_1)=\delta_1(T_1)+\delta_2(T_1).
\eeqlb
From (\ref{lem3-4}),(\ref{lem3-5}),(\ref{lem3-8}),(\ref{lem3-9}), (\ref{lem3-12}) and (\ref{lem3-14}), we obtain that
 \beqlb\label{lem3-15}
  d_{G}(T_1, T_2; s)&\leq& \delta(T_1)+M\int_0^{s}\frac{d_{G}(T_1, T_2; u)}{\sqrt{s-u}}\d u,
 \eeqlb
for every $s\in[0, 1]$. Furthermore, (\ref{lem3-15}) implies that
 \beqlb\label{lem3-16}
  d_{G}(T_1, T_2; s)&\leq& \delta(T_1)+M\int_0^{s}\frac{\delta(T_1)}{\sqrt{s-u}}\d u+M\int_0^s\frac{\d u}{\sqrt{s-u}}\int_0^{u}\frac{d_{G}(T_1, T_2; r)}{\sqrt{u-r}}\d r\nonumber
  \\&\leq&\delta(T_1)M+M\int_0^s d_{G}(T_1, T_2; r)\d r\int_r^{s}\frac{\d u}{\sqrt{(u-r)(s-u)}}.
 \eeqlb
Observe that
 \beqlb\label{lem3-21}
 \int_r^{s}\frac{\d u}{\sqrt{(u-r)(s-u)}}=\int_0^1\frac{\d y}{\sqrt{y-y^2}}=\pi.
 \eeqlb
(\ref{lem3-16}) yields that
 \beqlb\label{lem3-17}
  d_{G}(T_1, T_2; s)\leq\delta(T_1)M+ M\int_0^s d_{G}(T_1, T_2; r)\d r,
 \eeqlb
which and the Gronwall's inequality implies that
 \beqlb\label{lem3-18}
  d_{G}(T_1, T_2; s)\leq M\delta(T_1)\e^{M s},
 \eeqlb
for every $s\in[0, 1]$. Since (\ref{lem3-10}), (\ref{lem3-12}) and (\ref{lem3-14}) imply that
 $$\delta(T_1)\to 0,$$
as $T_1\to\infty$, (\ref{lem3-18}) indicates that for every $s\in[0, 1]$,
 \beqlb\label{lem3-19}
 d_{G}(T_1, T_2; s)\to 0,
 \eeqlb
for any $T_2>T_1$ and $T_1\to\infty$. From (\ref{lem3-19}) we derive that for every  $s\in[0, 1]$, $\{\int_{\R}G_T(x, s)\d x\}_T$ is a cauchy sequence. Therefore, for every $s\in[0, 1]$ there exists a function $G(s)$ such that
 \beqlb\label{lem3-20}
 \int_{\R}G_T(x, s)\d x\to G(s),
 \eeqlb
as $T\to\infty$. Note for any fixed $T$, $\int_{\R}G_T(x, s)\d x$ is non-negative and measurable on $s\in[0, 1]$. Furthermore, (\ref{lem7-1}) and the integrability of $\sigma(x)$ imply that for all $T>0$ $\int_{\R}G_T(x, s)\d x$ is bounded. Therefore $G(s)$ is measurable and
bounded as well.\qed

\begin{lem}\label{lem4} There exists a constant $M>0$ such that
 \beqlb\label{lem4-1}
  M\geq I_3(T, x, s),
 \eeqlb
for all $(x, s)\in\R\times[0, 1]$,  and as $T\to\infty$
 \beqlb\label{lem4-2}
 I_3(T, x, s)\to \gamma \int_0^{s}\frac{G(s-u)\d u}{\sqrt{2\pi u}},
 \eeqlb
for all $(x, s)\in\R\times[0, 1]$.
\end{lem}
{\bf Proof}\;\; From (\ref{s3-2-3}) it follows that
\beqlb\label{lem4-3}
I_3(T, x, s)=\gamma \int_0^{s}\frac{\d u}{\sqrt{2\pi u}}\int_{\R}\e^{-\frac{(x-y)^2}{2Tu}}\sigma(y)\bar{H}_T^2(y, s-u))\d y.
\eeqlb
(\ref{lem7-1}) and the integrability of $\sigma(x)$ yield that there is a constant $M>0$ such that
\beqlb
I_3(T, x, s)\leq M\int_0^{s}\frac{\d u}{\sqrt{2\pi u}}\leq M,
\eeqlb
for all $(x, s)\in\R\times[0, 1]$. In addition, (\ref{lem4-3}) implies that
\beqlb\label{lem4-4}
I_3(T, x, s)&=&\gamma \int_0^{s}\frac{\d u}{\sqrt{2\pi u}}\int_{\R}\sigma(y)\bar{H}_T^2(y, s-u)\d y\nonumber
\\&&-\gamma \int_0^{s}\frac{\d u}{\sqrt{2\pi u}}\int_{\R}(1-\e^{-\frac{(x-y)^2}{2Tu}})\sigma(y)\bar{H}_T^2(y, s-u)\d y.
\eeqlb
By using Lemma \ref{lem3}, it is easy to see that as $T\to\infty$,
\beqlb\label{lem4-5}
\gamma \int_0^{s}\frac{\d u}{\sqrt{2\pi u}}\int_{\R}\sigma(y)\bar{H}_T^2(y, s-u)\d y\to\gamma \int_0^{s}\frac{G(s-u)\d u}{\sqrt{2\pi u}},
\eeqlb
for all $(x, s)\in\R\times[0, 1]$. Furthermore, by using (\ref{lem7-1}), for all $(x, s)\in\R\times[0, 1]$
\beqlb\label{lem4-6}
&&\gamma \int_0^{s}\frac{\d u}{\sqrt{2\pi u}}\int_{\R}(1-\e^{-\frac{(x-y)^2}{2Tu}})\sigma(y)\bar{H}_T^2(y, s-u)\d y\nonumber
\\&&\qquad\leq M\int_0^{s}\frac{\d u}{\sqrt{2\pi u}}\int_{\R}(1-\e^{-\frac{(x-y)^2}{2Tu}})\sigma(y)\d y\to 0,
 \eeqlb
as $T\to\infty$. Applying (\ref{lem4-5}) and (\ref{lem4-6}) to  (\ref{lem4-4}), we arrive at (\ref{lem4-2}).\qed

\begin{lem}\label{lem5}
For any $(x, s)\in\R\times[0, 1]$, as $T\to\infty$,
 \beqlb\label{lem5-1}
 \bar{H}_T(x, s)\to \sum_{k=1}^n\int_0^s\frac{{\bf 1}_{[0, t_k]}(1-u)}{\sqrt{2\pi (s-u)}}\d u\int_{\R}\phi_k(y)\d y-\gamma \int_0^{s}\frac{G(s-u)\d u}{\sqrt{2\pi u}}.
 \eeqlb
Let $\Lambda(s)$ $(0\leq s\leq 1)$ denote the right hand side of (\ref{lem5-1}). Then $\Lambda(s)$ is the unique non-negative solution of the equation
 \beqlb\label{lem5-2}
 y(s)=\sum_{k=1}^n\int_0^s\frac{{\bf 1}_{[0, t_k]}(1-u)}{\sqrt{2\pi (s-u)}}\d u\int_{\R}\phi_k(y)\d y-K \int_0^{s}\frac{y^2(u)\d u}{\sqrt{2\pi (s-u)}}.\quad
 \eeqlb
\end{lem}
{\bf Proof}\; (\ref{lem5-1}) is an immediate conclusion from (\ref{s4-5}) and Lemma \ref{lem1}, Lemma \ref{lem2} and Lemma \ref{lem4}. $\Lambda(s)$ is non-negative and bounded because $\bar{H}_T(x, s)$ is non-negative and $0\leq \bar{H}_T(x, s)\leq M$ for all $(x, s)\in\R\times[0, 1]$ and $T>0$. To show $\Lambda(s)$ is a solution of (\ref{lem5-2}), we apply (\ref{lem5-1}) and the dominated convergence theorem to (\ref{lem3-1}), and then get that
 \beqlb\label{lem5-3}
 G(s)=\Lambda^2(s)\int_{\R}\sigma(y)\d y.
 \eeqlb
Substituting (\ref{lem5-3}) into (\ref{lem5-1}), we have that $\Lambda(s)$ satisfies the equation (\ref{lem5-2}). To prove the uniqueness, we suppose there is another non-negative solution $\Theta(s)$ of (\ref{lem5-2}). Let
 $$r(s):=|\Lambda(s)-\Theta(s)|, \qquad s\in[0, 1].$$
Note that the non-negative solution of (\ref{lem5-2}) should be bounded. There exists a constant $M>0$ such that
 $$|\Lambda^2(s)-\Theta^2(s)|\leq M r(s).$$
Therefore, from (\ref{lem5-2}) we get that
 $$
 r(s)\leq M\int_0^{s}\frac{r(s-u)\d u}{\sqrt u}.
 $$
Then
 $$
 r(s)\leq M\int_0^{s}\frac{r(u)\d u}{\sqrt{s-u}}\leq M\int_0^s r(t)\d t\int_t^s\frac{\d u}{\sqrt{(s-u)(u-t)}},
 $$
which and (\ref{lem3-21}) yield that
  \beqlb\label{lem5-4}
 r(s)\leq \pi M\int_0^s r(t)\d t.
 \eeqlb
By the Gronwall's inequality, $r(s)\equiv0$, and hence, $\Lambda(s)=\Theta(s)$ for all $s\in[0, 1]$.\qed

 \begin{rem}\label{s4-rem} From (\ref{lem5-2}), it is easy to see that $\Lambda(s)=0$ for all $s\in[0, 1-t_n)$.
 \end{rem}

Now we define
  \beqlb\label{s4-1}
 H_T(x):=T^{1/2}V_{\psi_T}(T^{1/2}x, T, 0),
\eeqlb
From (\ref{s3-9}) we obtain that
  \beqnn
 H_T(x)&=&T^{1/2}\int_0^{1}TL_{Ts}\Big(\psi_T(\cdot, Ts)\big(1-V_{\psi_T}(\cdot, T(1-s), Ts)\big)\Big)(T^{1/2}x)\d s\nonumber
 \\&&-\gamma T^{1/2} \int_0^{1}TL_{Ts}\Big(\sigma(\cdot)V_{\psi_T}^2(\cdot, T(1-s),
 Ts)\Big)(T^{1/2}x)\d s,
 \eeqnn
which combining with (\ref{s2-0-2}) further implies that
 \beqlb\label{s4-2}
 H_T(x)&=&T\int_0^{1}\d s\int_{\R}p_s(x-T^{-1/2}y)\psi_T(y, Ts)\d y\nonumber
 \\&&-T\int_0^1\d s\int_{\R}p_s(x-T^{-1/2}y)\psi_T(y, Ts)V_{\psi_T}(y, T(1-s), Ts)\d y\nonumber
 \\&&-\gamma  \int_0^{1}\d sT\int_{\R}p_s(x-T^{-1/2}y)\sigma(y)V_{\psi_T}^2(y, T(1-s),
 Ts)\d y.
 \eeqlb
Applying (\ref{s3-lem3-2}) to (\ref{s4-2}) leads to
 \beqlb\label{s4-3}
 H_T(x)&=&\int_0^{1}\d s\int_{\R}p_s(x-T^{-1/2}y)\psi(y, s)\d y\nonumber
 \\&&-T^{-1/2}\int_0^1\d s\int_{\R}p_s(x-T^{-1/2}y)\psi(y, s)\bar{H}_{T}(y, 1-s)\d y\nonumber
 \\&&-\gamma  \int_0^{1}\d s\int_{\R}p_s(x-T^{-1/2}y)\sigma(y)\big(\bar{H}_{T}(y, 1-s)\big)^2\d y.
 \eeqlb

\begin{lem}\label{lem6}Suppose $K=\gamma\int_{\R}\sigma(x)\d x<\infty$.
For any $x\in\R$, as $T\to\infty$,
\beqlb\label{lem6-1}
H_T(x)&\to& \sum_{k=1}^n\int_0^{t_k}p_s(x)\d s\int_{\R}\phi_k(y)\d y
-K\int_0^{1}p_s(x)\Lambda^2(1-s)\d s.\qquad
\eeqlb
Moreover,
\beqlb\label{lem6-5}
\lim_{T\to\infty}\int_{\R}H_T(x)\d x=\sum_{k=1}^nt_k\int_{\R}\phi_k(y)\d y
-K\int_0^{1}\Lambda^2(s)\d s.\qquad
\eeqlb
\end{lem}
{\bf Proof}\; By the dominated convergence theorem, it is easy to see that as $T\to\infty$,
 \beqlb\label{lem6-2}
 \int_0^{1}\d s\int_{\R}p_s(x-T^{-1/2}y)\psi(y, s)\d y\to\sum_{k=1}^n\int_0^{t_k}p_s(x)\d s\int_{\R}\phi_k(y)\d y,
 \eeqlb
for all $x\in R$. In addition, by (\ref{s4-5}) and the dominated convergence theorem again, we can readily get that as $T\to\infty$,
 \beqlb\label{lem6-3}
 T^{-1/2}\int_0^1\d s\int_{\R}p_s(x-T^{-1/2}y)\psi(y, s)\bar{H}_{\psi_T}(y, 1-s)\d y\to 0,
 \eeqlb
for all $x\in R$. Furthermore, by using Lemma \ref{lem5}, (\ref{lem7-1}) and the dominated convergence theorem, we find that as $T\to\infty$
 \beqlb\label{lem6-4}
&&\gamma  \int_0^{1}\d s\int_{\R}p_s(x-T^{-1/2}y)\sigma(y)\big(\bar{H}_{T}(y, 1-s)\big)^2\d y\nonumber
\\&&\qquad\to\gamma \int_0^{1}p_s(x)\d s\int_{\R}\sigma(y)\Lambda^2(1-s)\d y ,
 \eeqlb
for all $x\in R$. Applying (\ref{lem6-2})-(\ref{lem6-4}) to (\ref{s4-3}) we immediately obtain (\ref{lem6-1}). Moreover, using (\ref{s4-3}) again, we get that
\beqlb\label{lem6-6}
 \int_{\R}H_T(x)\d x&=&\int_{\R}\d x\int_0^{1}\d s\int_{\R}p_s(x-T^{-1/2}y)\psi(y, s)\d y\nonumber
 \\&&-T^{-1/2}\int_{\R}\d x\int_0^1\d s\int_{\R}p_s(x-T^{-1/2}y)\psi(y, s)\bar{H}_{T}(y, 1-s)\d y\nonumber
 \\&&-\gamma \int_{\R}\d x \int_0^{1}\d s\int_{\R}p_s(x-T^{-1/2}y)\sigma(y)\big(\bar{H}_{T}(y, 1-s)\big)^2\d y\nonumber
 \\&=&\sum_{k=1}^nt_k\int_{\R}\phi_k(y)\d y-T^{-1/2}\int_0^1\d s\int_{\R}\psi(y, s)\bar{H}_{T}(y, 1-s)\d y\nonumber
 \\&&-\gamma  \int_0^{1}\d s\int_{\R}\sigma(y)\big(\bar{H}_{T}(y, 1-s)\big)^2\d y.
\eeqlb
By (\ref{lem7-1}), it is easy to see that
 \beqlb\label{lem6-7}
 T^{-1/2}\int_0^1\d s\int_{\R}\psi(y, s)\bar{H}_{T}(y, 1-s)\d y\to 0,
 \eeqlb
as $T\to\infty$. Furthermore, combining Lemma \ref{lem4}, (\ref{lem7-1}) and the dominated convergence theorem yields that
 \beqlb\label{lem6-8}
 \gamma  \int_0^{1}\d s\int_{\R}\sigma(y)\big(\bar{H}_{T}(y, 1-s)\big)^2\d y\to\gamma  \int_0^{1}\d s\int_{\R}\sigma(y)\Lambda^2(1-s)\d y.
 \eeqlb
Combining (\ref{lem6-6}) with (\ref{lem6-7}) and (\ref{lem6-8}), we have (\ref{lem6-5}).\qed

Now we are at the place to give the proof of Theorem \ref{thm2}.

\noindent{\bf Proof of Theorem \ref{thm2}}\;\; For any given non-negative $\phi_1,\phi_2,\cdots,\phi_n\in\mathcal{S}(\R)$ and $0\leq t_1<t_2<\cdots<t_n\leq 1$, by using (\ref{thm2-1}), we obtain that
 \beqlb\label{thm2-2}
  &&\bfE\Big(\exp\Big\{-\sum_{k=1}^n \langle X(t_k), \phi_k\rangle\Big\}\Big)=\bfE\Big(\exp\Big\{-\sum_{k=1}^n \int_{\R}\phi_k(y)\d y \xi(t_k)\Big\}\Big)\nonumber
\\&&\qquad\qquad=\exp\bigg\{K\int_{1-t_n}^{1}\Lambda^2(s)\d s-\sum_{k=1}^nt_k\int_{\R}\phi_k(y)\d y\bigg\},
  \eeqlb
where $\Lambda(s)$ is exactly the unique nonnegative solution of the equation (\ref{lem5-2}). On the other hand, from Lemma \ref{s3-lem3} and (\ref{s4-1}), it follows that
 \beqlb\label{thm2-3}
 \bfE\Big(\exp\Big\{-\sum_{k=1}^n \langle X_T(t_k), \phi_k\rangle\Big\}\Big)=\exp\Big\{-\int_{\R}H_T(x)\d x\Big\}.
 \eeqlb
Applying Lemma \ref{lem6} and Remark \ref{s4-rem}, we obtain that
 \beqlb\label{thm2-4}
 &&\lim_{T\to\infty}\bfE\Big(\exp\Big\{-\sum_{k=1}^n \langle X_T(t_k), \phi_k\rangle\Big\}\Big)\nonumber
 \\&&\qquad=\exp\bigg\{K\int_{1-t_n}^{1}\Lambda^2(s)\d s-\sum_{k=1}^nt_k\int_{\R}\phi_k(y)\d y\bigg\}.
 \eeqlb
(\ref{thm2-2}) and (\ref{thm2-4}) imply that $X=\lambda\xi$ is the limit of $X_T$ in the sense of convergence of finite-dimensional distributions. The non-negativity of $\xi$ follows from the fact that $X_T$ is a non-negative measure for every $T>0$.\qed

Below we give the proof of Theorem \ref{thm1}.

\noindent{\bf Proof of Theorem 2.2}\:\:To prove the weak convergence of $\{X_T\}_{T\geq 1}$ in $C([0, 1], \mathcal{S}'(\R))$, it suffices to prove the convergence of finite-dimensional distributions plus the tightness of  $\{X_T\}_{T\geq 1}$ in $C([0, 1], \mathcal{S}'(\R))$. The former is proved by Theorem \ref{thm2}. To prove the latter, let $Y_T(t)=X_T(t)-t\lambda$. The tightness of $X_T$ in $C([0, 1], \mathcal{S}'(\R))$ is same as that of $Y_T$ in $C([0, 1], \mathcal{S}'(\R))$. Therefore, by using the
theorem of Mitoma \cite{M83}, it suffices to prove that $\{\langle Y_T, \phi\rangle\}_{T\geq 1}$ is tight in $C([0, 1],
 \R)$ for any given $\phi\in\mathcal{S}(\R)$. Thus, by the same arguments as those used in the proof of Theorem 2.1 in Li \cite{Li11c},
 we can readily get  the desired conclusion.
The details are omitted.\qed

\section{Properties of the solutions of integral equations}

In this section, we prove Theorem \ref{thm5}.

\medskip

\noindent{\bf Proof}. {\bf (1)}.
Let $\Lambda(s,\theta)$ be the unique nonnegative solution of the equation (\ref{thm5-1}) for $s\in[0, 1]$.
Therefore, $\Lambda(s,\theta)$ is continuous and differentiable. Furthermore, from the proof of Theorem \ref{thm2}, we see that $\Lambda(s,\theta)$ is the limit of
 \beqnn
 \bar{H}_T(x, s)&=&T^{1/2}V_{\psi_T}(x, Ts, T(1-s))\nonumber
 \\&=&T^{1/2}\Big[1-\bfE_x\Big(\exp\big\{-\int_0^{Ts}\big\langle N(u), \psi_T(\cdot, T(1-s)+u)\big\rangle\d u\big\}\Big)\Big]\nonumber
 \\&=&T^{1/2}\Big[1-\bfE_x\Big(\exp\big\{-\int_0^{s}\big\langle N(Tu), \phi(x){\bf 1}_{[0, 1]}(1-s+u)\big\rangle\d u\big\}\Big)\Big]\nonumber
 \\&=&T^{1/2}\Big[1-\bfE_x\Big(\exp\big\{-\int_0^{s}\big\langle N(Tu), \phi(x)\big\rangle\d u\big\}\Big)\Big],
 \eeqnn
where $\phi\geq 0$ satisfies that $\int_{\R}\phi(x)\d x=\theta$.
It is easy to see that $\bar{H}_T(x, s)$ is non-decreasing on $s$. Therefore $\Lambda(s,\theta)$ is non-decreasing as well. Differentiating both sides of (\ref{thm5-1}) on $s$ leads to that for any $s\in(0, 1)$
\beqlb\label{prop-3}
\Lambda'(s,\theta)=\sqrt{\frac{1}{2\pi s}}\theta-2 K \int_0^{s}\frac{\Lambda(s-u,\theta)\Lambda'(s-u,\theta)\d u}{\sqrt{2\pi u}}\geq0,\quad
\eeqlb
which further implies that
\beqlb\label{prop-4}
\sqrt{2\pi s}\Lambda'(s,\theta)&=&\theta-2K\sqrt{s}\int_0^{s}\frac{\Lambda(s-u,\theta)\Lambda'(s-u,\theta)\d u}{\sqrt{u}}\leq \theta-K\Lambda^2(s,\theta).
 \eeqlb
Consequently, we have that
 \beqnn
 \Lambda(s,\theta)\leq \sqrt{\theta/K}.
 \eeqnn

{\bf (2)}. Now, we observe the following nonlinear integral equation.
 \beqlb\label{prop-3-0}
 h(s)=\sqrt{\frac{2s}{\pi}}-K \int_0^{s}\frac{h^2(u)\d u}{\sqrt{2\pi (s-u)}},\qquad s\geq 0.
 \eeqlb
We can readily verify that for any $m>0$,
 \beqlb\label{prop-3-01}
 h(s):=\Lambda(s/\theta,\theta)/\sqrt{\theta}\leq1/\sqrt{K}
 \eeqlb
  is a non-negative and bounded solution of (\ref{prop-3-0}) on $[0, m]$ for a given $\theta\geq m$. On the contrary, if $h(s)$ is a solution of (\ref{prop-3-0}) on $[0, m]$, then for any $\theta\leq m$,
 \beqlb\label{prop-3-1}
 \Lambda(s,\theta)=\sqrt{\theta}h(\theta s), \qquad s\in[0, 1],
 \eeqlb
is a non-negative and bounded solution of (\ref{thm5-1}) for $s\in[0, 1]$. Due to the uniqueness of non-negative and bounded solutions of (\ref{thm5-1}) for $s\in[0, 1]$, we know (\ref{prop-3-0}) has an unique non-negative and bounded solution on $[0, m]$ for any $m>0$, and hence on $[0,\infty)$. For convenience, we denote the unique solution by $h(s)$. Thanks to the monotonicity of $\Lambda(s,\theta)$ on $s\in[0, 1]$ and the definition of (\ref{prop-3-01}), we can readily see that $h(s)$ is non-decreasing in $[0. m]$ for any $m>0$. Therefore, $h(s)$ is non-decreasing on $[0,\infty)$.

Based on the aforementioned facts, we can extend $\Lambda(s,\theta)$ to the positive half-line by (\ref{prop-3-1}) and denote the extension by $\Lambda(s,\theta)$ as well. Obviously, $h(s)=\Lambda(s, 1)$. It is easy to check that $\Lambda(s,\theta)$ is the unique non-decreasing non-negative solution of (\ref{thm5-1}) for all $s\geq0$. Furthermore by using (\ref{prop-3-01}) we obtain that for any $s>0$ and $\theta>0$
 $$\Lambda(s,\theta)=\sqrt{\theta}h(\theta s)=\frac{\sqrt{\theta}}{\sqrt{\theta s}}\Lambda(\frac{\theta s}{\theta s}, \theta s)=\frac{1}{\sqrt{s}}\Lambda(1,\theta s).$$
Therefore, the desired conclusion (\ref{thm5-2}) holds.  Let
  $$l=\lim_{s\to\infty}h(s)=\lim_{s\to\infty}\Lambda(s, 1).$$
Then from (\ref{prop-3-0}) we derive that
 \beqnn
 l\geq \sqrt{\frac{2s}{\pi}}-Kl^2 \int_0^{s}\frac{\d u}{\sqrt{2\pi (s-u)}}=\sqrt{\frac{2s}{\pi}}(1-Kl^2),
 \eeqnn
for any $s>0$. Letting $s\to\infty$ leads to
 \beqlb\label{prop-3-2}
 Kl^2\geq 1.
 \eeqlb
On the other hand,  from Theorem \ref{thm1} it follows that
 \beqlb\label{prop-3-3}
1\geq\bfE(\e^{-\theta\xi(1)}) &=&\exp\bigg\{K\int_0^{1}\Lambda^2(s,\theta)\d s-\theta \bigg\}=\exp\bigg\{\theta\Big(K\int_0^1 \Lambda^2(\theta s, 1)\d s-1\Big)\bigg\}\nonumber
\\&=&\exp\bigg\{\theta\Big(K\int_0^\theta \Lambda^2(s, 1)\d s/\theta-1\Big)\bigg\},
 \eeqlb
which indicates that
 \beqnn
\lim_{\theta\to\infty}K\int_0^\theta \Lambda^2(s, 1)\d s/\theta-1=Kl^2-1\leq 0.
 \eeqnn
Therefore, we have that
 \beqlb\label{prop-3-4}
 \lim_{s\to\infty}\Lambda(s, 1)=l=1/\sqrt{K}.
 \eeqlb
The proof of Theorem \ref{thm5} is complete.\qed

 \section{Properties of the limiting process $\xi$}

By the same discussion as those in the proof of Talarczyk \cite[Theorem 2.4]{T08}, we can readily get  from Theorem \ref{thm1} that the process $\xi$ is infinitely divisible, self-similar with index $1$, non-decreasing and nonnegative, and has continuous paths, i.e., Proposition \ref{prop4} is right. Below, we prove Proposition \ref{prop1}.

\noindent{\bf Proof}.\;Recall that $\xi(1)$ is non-negative and infinitely divisible. From \cite[P.385]{S99}, we know that there exist a non-negative constant $\gamma_0$ and a measure $\nu$ on $(0,\infty)$ with $\int_0^\infty (1\wedge x)\nu(\d x)<\infty$ such that
 \beqlb\label{prop-3-5}
\bfE(\e^{-\theta\xi(1)})=\exp\bigg\{-\theta\Big(\gamma_0+\int_0^\infty\frac{1-\e^{-\theta x}}{\theta}\nu(\d x)\Big)\bigg\},
 \eeqlb
which and (\ref{prop-3-3}) lead to
 \beqlb\label{prop-3-6}
 K\int_0^\theta \Lambda^2(s,1)\d s/\theta-1=-\gamma_0-\int_0^\infty\frac{1-\e^{-\theta x}}{\theta}\nu(\d x).
 \eeqlb
 Note that $\int_0^\infty\frac{1-\e^{-\theta x}}{\theta}\nu(\d x)\to 0$ as $\theta\to\infty$. From (\ref{prop-3-4}) and (\ref{prop-3-6}) we obtain that
  \beqlb\label{prop-3-7}
  \gamma_0=0.
  \eeqlb
 Substituting  (\ref{prop-3-7}) into (\ref{prop-3-5}) yields
 \beqnn
\bfE(\e^{-\theta\xi(1)})=\exp\bigg\{-\int_0^\infty(1-\e^{-\theta x})\nu(\d x)\bigg\}.
 \eeqnn
Therefore, (\ref{prop-1}) follows from the fact  $\xi(t)\stackrel{d}{=}t\xi(1)$.

 Now, substituting  (\ref{prop-3-7}) into (\ref{prop-3-6}), we get that
  \beqlb\label{prop-3-8}
    K\int_0^\theta \Lambda^2(s,1)\d s=\theta-\int_0^\infty(1-\e^{-\theta x})\nu(\d x).
 \eeqlb
 Differentiating both sides of (\ref{prop-3-8}) leads to
  \beqlb\label{prop-3-9}
  K\Lambda^2(\theta,1)=1-\int_0^\infty x\e^{-\theta x}\nu(\d x),
  \eeqlb
for any $\theta>0$. Combining (\ref{thm5-2}) with (\ref{prop-3-9}) indicates (\ref{thm5-4}). 

In addition, substituting (\ref{prop-3-9}) into (\ref{prop-3-0}) leads to
  \beqlb\label{prop-3-13}
  \Lambda(t, 1)&=&\int_0^t\frac{\d \theta}{\sqrt{2\pi (t-\theta)}}\Big(\int_0^\infty x\e^{-\theta x}\nu(\d x)\Big)=\int_0^\infty \Big[x\e^{-t x}\int_0^t\frac{\e^{\theta x}}{\sqrt{2\pi \theta}}\d \theta\Big]\nu(\d x)\nonumber
  \\&=&\int_0^\infty \frac{1}{\sqrt{2\pi t}}Q(tx)\nu(\d x),
  \eeqlb
where for any $w>0$, $$Q(w)=\frac{\sqrt w}{\e^w}\int_0^{w}\frac{\e^{y}}{\sqrt y}\d y.$$
Therefore,  combining (\ref{thm5-2}) and (\ref{prop-3-13}) leads to (\ref{thm5-3}). 

Furthermore, it is easy to see that $Q(w)$ is continuous on $(0,\infty)$ with $\lim_{w\to 0}Q(w)=0$, and by the L'H\^{o}pital's law,
 $$\lim_{w\to\infty}Q(w)=1.$$
Therefore, there exists a constant $M>0$ such that $Q(w)\leq M$ for all $w>0$. Applying this fact to (\ref{prop-3-13}), we get that
 \beqlb
 \sqrt{2\pi t}\Lambda(t,1)\leq M\int_0^\infty \nu(\d x),
 \eeqlb
for all $t>0$. Then letting $t\to\infty$ implies that
 \beqlb\label{prop-3-14}
 \nu((0,\infty))=\int_0^\infty \nu(\d s)=\infty.
 \eeqlb
In addition, from the equation (\ref{thm5-1}), it is easy to see that $\Lambda(0, 1)=0$. Letting $\theta\to 0$, (\ref{prop-3-9}) implies that
 \beqlb\label{prop-3-12}
 1=\int_0^\infty x\nu(\d x).
 \eeqlb
Moreover, (\ref{prop-3-0}) implies that $h(s):=\Lambda(s, 1)$ is  differentiable, and for all $s>0$
 \beqlb\label{prop-3-10}
 h'(s)=\frac{1}{\sqrt{2\pi s}}-K\int_0^s\frac{2h(s-u)h'(s-u)}{\sqrt{2\pi u}}\d u.
 \eeqlb
Applying (\ref{prop-3-10}) to (\ref{prop-3-0}) yields that for all $s>0$
 \beqnn
 h(s)h'(s)&=& \frac{1}{\pi}-K\sqrt{\frac{2s}{\pi}}\int_0^s\frac{2h(s-u)h'(s-u)}{\sqrt{2\pi u}}\d u-\frac{K}{\sqrt{2\pi s}}\int_0^s\frac{h^2(s-u)}{\sqrt{2\pi u}}\d u\nonumber
 \\&&+K^2\int_0^s\frac{2h(s-u)h'(s-u)}{\sqrt{2\pi u}}\d u\int_0^s\frac{h^2(s-u)}{\sqrt{2\pi u}}\d u.
 \eeqnn
Observe that by the L'h\^{o}pital Law,
 $$\lim_{s\to 0}\frac{1}{\sqrt{2\pi s}}\int_0^s\frac{h^2(s-u)}{\sqrt{2\pi u}}\d u=\lim_{s\to 0}\sqrt{\frac{2s}{\pi}}\int_0^s\frac{2h(s-u)h'(s-u)}{\sqrt{2\pi u}}\d u=0.$$
It is easy to see that
 \beqlb\label{prop-3-11}
 \lim_{s\to 0} h(s)h'(s)=\frac{1}{\pi}.
 \eeqlb
Now, we differentiate both sides of (\ref{prop-3-9}) and get that
 \beqnn
  2Kh(\theta)h'(\theta)=\int_0^\infty x^2\e^{-\theta x}\nu(\d x).
 \eeqnn
Letting $\theta\to0$ and using (\ref{prop-3-11}), we obtain that
 \beqnn
 \int_0^\infty x^2\nu(\d x)=2K/\pi,
 \eeqnn
which  and (\ref{prop-3-14}), and (\ref{prop-3-12}) complete the proof of Proposition \ref{prop1}.\qed

At last, we devote to proving that $X_T$ converges in finite-dimensional distributions to the $0$-measure under the condition of $\int_{\R}\sigma(x)\d x=\infty$. This is equivalent to proving that $X_T(1)$ converges in distributions as $T\to\infty$ to the measure concentrated on the $0$-measure.

{\bf Proof of Theorem \ref{thm4}}.\;\; Let $B_n=\{|x|\leq n\}$ and $\sigma_n(x)=\sigma(x){\bf 1}_{B_n}(x)$. Then $$K_n:=\gamma\int_{\R}\sigma_n(x)\d x<\infty.$$
For any non-negative $\phi\in\mathcal{S}(\R)$, Let $v_n(x, t)\in[0, 1]$ be the solution of
\beqnn
 \frac{\partial v(x, t)}{\partial t}&=&
 \frac{1}{2}\Delta v(x, t)-\gamma\sigma_n(x)v^2(x, t)+\frac{1}{T}\phi(x){\bf 1}_{[0, 1]}(1-\frac{t}{T})\big(1-v(x, t)\big),
 \eeqnn
with $v(x,0)=0$, and   $v(x, t)\in[0, 1]$ be also the solution of this equation but $\sigma_n$ replaced by $\sigma$.
Since $0\leq\sigma_n\leq\sigma$, by the maximum principle, we have that
 $$
 0\leq v(x, t)\leq v_n(x, t)\leq 1.
 $$
Therefore, from Remark \ref{s3-rem}, it follows that for any $n\geq 1$,
 \beqlb\label{s6-1}
 \bfE\Big(\exp\Big\{-\langle X_T(1), \phi\rangle\Big\}\Big)\geq\exp\bigg\{-\int_{\R} v_n(x, T)\d x\bigg\}.
 \eeqlb
Reviewing the proof of Theorem \ref{thm2}, we have that
 \beqlb\label{s6-2}
 \lim_{T\to\infty}\exp\bigg\{-\int_{\R} v_n(x, T)\d x\bigg\}=\exp\bigg\{K_n\int_0^{1}\Lambda^2_n(s)\d s-\|\phi\|_1\bigg\},
 \eeqlb
where the function $\Lambda_n(s)$ is the unique non-negative solution of the equation
 \beqlb\label{s6-3}
 \Lambda(s)=\sqrt{\frac{2s}{\pi}}\|\phi\|_1-K_n\int_0^{s}\frac{\Lambda^2(u)\d u}{\sqrt{2\pi (s-u)}},\quad
 \eeqlb
and $\|\phi\|_1=\int_{\R}\phi(y)\d y$. Let $\bar{\Lambda}_n(s)=K_n\Lambda_n(s)$. Then (\ref{s6-3}) implies that $\bar{\Lambda}_n(s)$ is the unique non-negative solution of the equation
 \beqlb\label{s6-4}
 \Lambda(s)=\sqrt{\frac{2s}{\pi}}K_n\|\phi\|_1-\int_0^{s}\frac{\Lambda^2(u)\d u}{\sqrt{2\pi (s-u)}}.
 \eeqlb
Observing the equation (\ref{thm5-1}), we get that
 \beqnn
 \bar{\Lambda}_n(s)=\Lambda_1(s, K_n\|\phi\|_1),
 \eeqnn
where $\Lambda_1(s, \theta)$ is the unique non-negative solution of (\ref{thm5-1}) with $K=1$. Then by same argument used in the proof of Proposition \ref{prop1}, we can readily get that
 \beqlb\label{s6-5}
 K_n\int_0^{1}\Lambda^2_n(s)\d s&=&\frac{1}{K_n}\int_0^{1}\bar{\Lambda}^2_n(s)\d s\nonumber
 \\&=&\frac{1}{K_n}\int_0^{1}\Lambda_1^2(s, K_n\|\phi\|_1)\d s=\frac{1}{K_n}\int_0^{K_n\|\phi\|_1}\Lambda_1^2(s, 1)\d s.
 \eeqlb
By (\ref{prop-3-4}),
 \beqlb\label{s6-6}
 \lim_{s\to\infty} \Lambda_1(s, 1)=1.
 \eeqlb
The assumption $\int_{\R}\sigma(x)\d x=\infty$ implies that $K_n\to\infty$ as $n\to\infty$. From (\ref{s6-5}) and (\ref{s6-6}), it follows that
 \beqlb\label{s6-7}
 \lim_{n\to\infty}K_n\int_0^{1}\Lambda^2_n(s)\d s=\|\phi\|_1.
 \eeqlb
Combining (\ref{s6-1}) with (\ref{s6-2}) and letting $n\to\infty$, we get that
  \beqnn
 \lim_{T\to\infty}\bfE\Big(\exp\Big\{-\langle X_T(1), \phi\rangle\Big\}\Big)&\geq &\lim_{n\to\infty}\lim_{T\to\infty}\exp\bigg\{-\int_{\R} v_n(x, T)\d x\bigg\}
 \\&=&\lim_{n\to\infty}\exp\bigg\{K_n\int_0^{1}\Lambda^2_n(s)\d s-\|\phi\|_1\bigg\}=1,
 \eeqnn
where (\ref{s6-7}) is used at the last equality. The proof of Theorem \ref{thm4} is complete.\qed 

\end{document}